\documentclass[final,leqno,onefignum,onetabnum]{siamltex1213}
\usepackage{algorithm}
\usepackage{algorithmic}
\usepackage{amsfonts}
\usepackage{amsmath}
\usepackage{array}
\usepackage{hyperref}
\hypersetup{colorlinks,linkcolor=black,urlcolor=gray,citecolor=black}
\usepackage{tikz}
\usetikzlibrary{shapes.geometric, arrows, calc}

\title{A Parallel Orbital-updating Based Optimization Method for Electronic Structure Calculations \thanks{This work was supported by the National Science Foundation of China under grant 9133202, the Funds for Creative Research Groups of
China under grant 11321061,  the National Basic Research Program
of China under grant 2011CB309703, and  and the National Center for Mathematics and Interdisciplinary Sciences of the Chinese Academy of Sciences
.}}

\begin{document}
\author{
Xiaoying Dai\thanks{LSEC, Institute of Computational Mathematics
and Scientific/Engineering Computing, Academy of Mathematics and
Systems Science, Chinese Academy of Sciences, Beijing 100190, China
(daixy@lsec.cc.ac.cn).}, \and  Zhuang Liu\thanks{LSEC, Institute of Computational Mathematics
and Scientific/Engineering Computing, Academy of Mathematics and
Systems Science, Chinese Academy of Sciences, Beijing 100190, China
(liuzhuang@lsec.cc.ac.cn).},\and Xin Zhang\thanks{LSEC, Institute of Computational Mathematics
and Scientific/Engineering Computing, Academy of Mathematics and
Systems Science, Chinese Academy of Sciences, Beijing 100190, China
(xzhang@lsec.cc.ac.cn).}, \and Aihui Zhou\thanks{LSEC, Institute of Computational Mathematics
and Scientific/Engineering Computing, Academy of Mathematics and
Systems Science, Chinese Academy of Sciences, Beijing 100190, China
(azhou@lsec.cc.ac.cn).}
}

\maketitle


\begin{abstract}
In this paper, we propose a parallel optimization method
for electronic structure calculations based on a single orbital-updating approximation. It is shown by our numerical experiments that
the method is efficient and reliable  for atomic and molecular systems of large scale over supercomputers.
\end{abstract}



\pagestyle{myheadings}
\thispagestyle{plain}

\section{Introduction}


Kohn-Sham equations are nonlinear eigenvalue problems and typical models used in
electronic structure calculations. Due to the nonlinearity and requirement
of a cluster of eigenpairs, it is desired and significant to design efficient and
in particular supercomputer friendly numerical methods for solving Kohn-Sham equations.
Most recently, a so-called parallel orbital-updating algorithm has been proposed in
\cite{DGZZ} to Kohn-Sham equations and been demonstrated to be a quite potential
approach for electronic structure calculations over supercomputers.

Note that in the context of solving nonlinear eigenvalue problems, the self consistent
field (SCF) iteration is most widely used. However, the convergence of SCF is not
guaranteed, especially for large scale systems with small band gaps, its performance is
unpredictable \cite{ZZWZ}. As a result, optimization algorithms have been developed for
the Kohn-Sham direct energy minimization models, which are different from solving the nonlinear
eigenvalue models, including trust region type methods
\cite{FMM, FMM2, SH, TOYJSH},  Newton-type methods \cite{EAS, LWWUY},
 conjugate gradient methods \cite{PTLJ, SCPB, WVHN}, and gradient-type methods \cite{UWYKL,WY, ZZWZ}
 (see, also, \cite{Sa,YaMeWa} and references cited therein).
We particularly mention that the gradient type method proposed in \cite{ZZWZ} constructs new trail
points along the gradient on the Stiefel manifold and has been shown to be efficient and robust.
However,
 to develop a parallel optimization algorithm is still significant and challenging  for electronic structure calculations, especially
for systems of large scale over supercomputers.


We observe that optimization algorithms often have better convergence properties:
as long as the ground state is not orthogonal to the starting point, the ground state
is reachable from any starting point by following a path of decreasing energy \cite{TePaAl}.
In this paper, we shall apply a similar idea of \cite{DGZZ} to direct energy minimization models
and propose a parallel approach to electronic structure
calculations based on a single orbital-updating approximation.



We understand that the parallel orbital-updating algorithm proposed in \cite{DGZZ} reduces
the large scale eigenvalue problem to the solution of some independent
source problems and some eigenvalue problem of small systems, and the source
problems are independent each other and can be calculated in parallel intrinsically.
While in this paper, instead of solving the large global energy  minimization problem, we simply
minimize the total energy functional along each orbital direction independently and then orthonormalize
the orbital minimizers obtained, which is a two-level parallelism similar to \cite{DGZZ}. We should
mention that in practice, it is not easy and unnecessary to find the minimum for each orbital
in each iteration. Instead, we solve the $N$ optimization subproblems approximately or inexactly.

The rest of this paper is organized as follows: in Section 2, we provide a brief introduction
to the Kohn-Sham density functional theory (DFT) model and  its associated Stiefel manifold.
We then propose our parallel orbital-updating based optimization algorithms for electronic structure
calculations in Section 3. We present numerical results in Section 4 to show the accuracy and
efficiency of our algorithms. Finally, we give some concluding remarks in Section 5.

\section{Preliminaries}
\subsection{Kohn-Sham DFT model}
Let $\mathcal{U} = (u_1, \dots, u_N)$, with $u_i\in H^1(\mathbb{R}^3)$, the Kohn-Sham DFT model
solves the following constrained optimization problem
\begin{equation}\label{emin}
\begin{split}
&\inf \ \ \ E(\mathcal{U}) \\
s.t.\ \  &\int_{\mathbb{R}^3} u_i(r)u_j(r) dr = \delta_{ij}, \ \ 1\leq i,j\leq N,
\end{split}
\end{equation}
where the Kohn-Sham total energy $E(\mathcal{U})$ is defined as
\begin{eqnarray}\label{energy}
E(\mathcal{U})&=&\frac{1}{2} \int_{\mathbb{R}^3} \sum_{i=1}^N|\nabla u_i(r)|^2 dr +
\frac{1}{2}\int_{\mathbb{R}^3}\int_{\mathbb{R}^3}\frac{\rho(r)\rho(r')}{|r-r'|}drdr' \nonumber\\
&&+\int_{\mathbb{R}^3} \sum_{i=1}^N u_i(r)V_{ext}(r)u_i(r)dr+ \int_{\mathbb{R}^3}
\varepsilon_{xc}(\rho(r))\rho(r)dr.
\end{eqnarray}
Here $\rho(r)=\sum \limits_{i=1}^N|u_i(r)|^2$ is the electronic density,
$V_{ext}(r)$ is the external potential generated by the nuclei: for full potential calculations,
$V_{ext}(r) = -\sum\limits_{I = 1}^{M}\frac{Z_I}{|r - R_I|}$, $Z_I$ and $R_I$ are the nuclei charge
and position of the $I$-th nuclei respectively; while for pseudo potential
approximations, $V_{ext}(r) = \sum\limits_{I = 1}^{M} V^I_{loc}(r) + V^I_{nloc}(r)$,
$V^I_{loc}(r)$ is the local part of the pseudo potential for nuclei $I$, $V^I_{nloc}(r)$
is the nonlocal part, which usually has the following form
\begin{equation*}
V^I_{nloc} u (r)= \sum_l\int_{\mathbb{R}^3}\xi_l^{I}(r') u(r') dr'\xi_l^{I}(r),
\end{equation*}
with $\xi_l^{I}\in L^2(\mathbb{R}^3)$ \cite{Ma}. The $\varepsilon_{xc}(\rho(r))$ in the forth
term is the exchange-correlation functional, describing the many-body effects of exchange and
correlation \cite{Ma}, which is not known explicitly, and some approximation (such as LDA, GGA)
has to be used. The functions $\{u_i, i = 1, \dots, u_N\}$ are call Kohn-Sham orbitals.

Denote $E_{u_i}$ the derivative of $E(\mathcal{U})$ to the $i$-th orbital, it is easy to see
that
\begin{equation*}
E_{u_i} = H(\rho) u_i,
\end{equation*}
where
\begin{equation*}
H(\rho) = -\frac{1}{2}\Delta + V_{ext} + \int_{\mathbb{R}^3}\frac{\rho(r')}{|r-r'|}dr' +
v_{xc}(\rho),
\end{equation*}
here $v_{xc}(\rho) = \frac{\text{d}(\rho\varepsilon_{xc}(\rho))}{\text{d}\rho}$ is the exchange
correlation potential, $H(\rho)$ is the Kohn-Sham Hamiltonian operator.

\subsection{Stiefel manifold and gradient}
The feasible set of the constraint problem \eqref{emin} is a Stiefel manifold, which
is defined as
\begin{equation*}
\mathcal{M}_{V}^N = \{\mathcal{U} = (u_i)_{i=1}^N|u_i\in V,
\langle\mathcal{U}^T\mathcal{U}\rangle = I_N\},
\end{equation*}
where $V = H^1(\mathbb{R}^3)$, and for $\Psi = (\psi_1,\dots,\psi_N)\in V^N$, $\Phi =
(\phi_1,\dots,\phi_N)\in V^N$, the inner product matrix is defined as $\langle\Psi^T\Phi\rangle
:= (\langle\psi_i,\phi_j\rangle)_{ij=1}^{N}\in \mathbb{R}^{N\times N}$, $\langle\psi_i,\phi_j\rangle
= \int_{\mathbb{R}^3} \psi_i(r)\phi_j(r) dr$
is the usual $L_2$ inner product. By computing the derivative of the constraint condition
\begin{equation*}
\langle\mathcal{U}^T\mathcal{U}\rangle = I_N,
\end{equation*}
we get the tangent space of $\mathcal{M}_{V}^N$ at $\mathcal{U}$
\begin{equation*}
\mathcal{T}_\mathcal{U}\mathcal{M}_{V}^N = \{ \mathcal{B}\in V^{N}: \langle \mathcal{U}^T\mathcal{B}
\rangle = -\langle \mathcal{B}^T\mathcal{U}\rangle \in \mathbb{R}^{N\times N}\}.
\end{equation*}

We verify from \eqref{energy} that $E(\mathcal{U}) =
E(\mathcal{U}P)$ for any $N\times N$ orthogonal matrix $P$. To get rid of the non uniqueness,
we consider the problem on the Grassmann manifold, which is the quotient of the Stiefel
manifold and defined as
\begin{equation*}
\mathcal{G}_{V}^N = \mathcal{M}_{V}^N/\sim,
\end{equation*}
here $\mathcal{U}\sim \tilde{\mathcal{U}}$ if there is an orthogonal matrix $P$ such that
$\tilde{\mathcal{U}} = \mathcal{U}P$, we use $[\mathcal{U}]$ to denote the equivalence class.
The tangent space of the Grassmann manifold $\mathcal{G}_{V}^N$ at $[\mathcal{U}]$ is \cite{SRNB}
\begin{equation*}
\mathcal{T}_{[\mathcal{U}]}\mathcal{G}_{V}^N = \{ \mathcal{F}\in V^{N}: \langle
\mathcal{F}^T\mathcal{U}\rangle = 0 \in \mathbb{R}^{N\times N}\}.
\end{equation*}

By the first order necessary condition, we see that the derivative $E'(\mathcal{U})$ of the energy
$E(\mathcal{U})$ at the minimizer $[\mathcal{U}]\in\mathcal{G}_{V}^N$ satisfies \cite{SRNB}
\begin{equation*}
\text{tr}\langle (E'(\mathcal{U}))^T\delta\mathcal{U}\rangle = \text{tr}\langle (H(\rho)\mathcal{U})^T
\delta\mathcal{U}\rangle = 0, \ \ \ \forall\ \delta\mathcal{U}
\in \mathcal{T}_{[\mathcal{U}]}\mathcal{G}_{V}^N,
\end{equation*}
i.e. the derivative $E'(\mathcal{U})$ vanishes on the tangent space
$\mathcal{T}_{[\mathcal{U}]}\mathcal{G}_{V}^N$ of the Grassmann manifold,
we then get the gradient $\nabla E(\mathcal{U})$ of the energy on the Stiefel manifold
\begin{equation}\label{gra}
\nabla E(\mathcal{U}) = H(\rho)\mathcal{U} - \mathcal{U}\Sigma = 0,
\end{equation}
where $\Sigma = \langle (H(\rho)\mathcal{U})^T\mathcal{U} \rangle$ is symmetric since the Hamiltonian
$H(\rho)$ is a symmetric operator. Suppose orthogonal matrix $P$ diagonalize $\Sigma$, i.e.
$P^T\Sigma P = \Lambda$, where $\Lambda$ is a diagonal matrix, let $\mathcal{U} = \mathcal{U}P$,
substitute into \eqref{gra}, notice that $\rho = \sum\limits_{i=1}^N|u_i|^2$ is invariant under the
orthogonal transformation, we get the well known Kohn-Sham equation
\begin{equation}
H(\rho)\mathcal{U} = \mathcal{U}\Lambda, \ \ s.t. \ \ \langle \mathcal{U}^T \mathcal{U} \rangle = I_N.
\end{equation}

\subsection{Discretization}
The continuous Kohn-Sham DFT model can be discretized by either real space approaches including
finite difference, finite element ect., or by planewave basis. We focus on the real space discretization
in this paper. For simplicity, we still use the notations as in the previous part for the
discretized problem. Note that this does not affect the algorithm we propose in this paper.

As will be shown in section 3, since the mesh and hence the finite dimensional space is changing, we
denote the discretized finite dimensional space by $V^{(n)}$ in outer iteration $n$. For the inner loop,
the solution $\mathcal{W}^{(l)}$ is in the space $V^{(n)}$.
\section{Parallel orbital-updating based optimization algorithm}

\subsection{General framework}
Motivated by \cite{DGZZ}, we propose a parallel orbital-updating based optimization
algorithm, whose general framework is shown as \textbf{Algorithm} \ref{Algorithm B}.

\begin{algorithm}
\caption{}
\label{Algorithm B}
\begin{algorithmic}[1]
\STATE Given initial space $V^{(0)}$, initial data $\mathcal{U}^{(0)}, \ s.t. \
       \langle(\mathcal{U}^{(0)})^T \mathcal{U}^{(0)}\rangle = I_N$, let $n = 0$;
\STATE Apply an adaptive approach to $ \mathcal{U}^{(n)}$ and get $V^{(n+1)}$;
\STATE $\mathcal{W}^{(0)} = \mathcal{U}^{(n)}$, $l = 0$;
\STATE For $i = 1,\dots, N$, apply an optimization approach to get $\widetilde{w}_i^{(l+1)}$ in parallel,
       such that
       \begin{equation*}
        E(\mathcal{W}^{(l)}_i\diamond {\tilde w}^{(l+1)}_i) < E(\mathcal {W}^{(l)}),
       \end{equation*}
       where $\mathcal{W}^{(l)}_i\diamond v = (w_1^{(l)},\dots,w_{i-1}^{(l)},v,w_{i+1}^{(l)},\dots,w_N^{(l)})$;
\STATE $\mathcal{W}^{(l+1)} = \text{Orth}(\widetilde{\mathcal{W}}^{(l+1)}), \ l = l + 1$;
\STATE Convergence check for the inner iteration: if not converged, goto step 4;
\STATE $\mathcal{U}^{(n+1)} = \mathcal{W}^{(l)}$, $n = n + 1$;
\STATE Convergence check for the outer iteration: if not converged, goto step 2;
       else, stop.
\end{algorithmic}
\end{algorithm}

\tikzstyle{basic} = [rectangle, minimum width=15pt, minimum height=10pt,text centered,
                     draw=black]
\tikzstyle{arrow} = [thick,->,>=stealth]
\tikzstyle{larrow} = [draw, -latex',thick]

\begin{figure}
\caption{Flowchart for \textbf{Algorithm} \ref{Algorithm B}}
\begin{center}
\begin{tikzpicture}[node distance=35pt] \label{Flowchart}
\node (start) [basic] {Given initial data $\mathcal{U}^{(0)}, \ s.t. \ \langle(\mathcal{U}^{(0)})^T
                       \mathcal{U}^{(0)}\rangle = I_N$, initial mesh $\mathcal{T}^{(0)}$, and let $n=0$};
\node (step1) [basic, below of=start] {Apply an adaptive approach to $ \mathcal{U}^{(n)}$ and get
                       $\mathcal{T}^{(n+1)}$};

\node (step2) [basic, below of=step1] {$\mathcal{W}^{(0)} = \mathcal{U}^{(n)}$, $l = 0$};

\node (step3) [basic, below of=step2] {Apply an optimization approach in parallel,
                                      and get ${\widetilde{\mathcal{W}}}^{(l+1)}$};
\node (step4) [basic, below of=step3] {Orthogonalization: $\mathcal{W}^{(l+1)} = \text{Orth}
                                      (\widetilde{\mathcal{W}}^{(l+1)}), \ l = l + 1$};
\node (step5) [basic, below of=step4] {Convergence check for the inner iteration};
\node (step6) [basic, below of=step5] {$\mathcal{U}^{(n+1)} = \mathcal{W}^{(l)}$, $n = n + 1$};
\node (step7) [basic, below of=step6] {Convergence check for the outer iteration};
\node (step8) [basic, below of=step7] {Stop};

\draw [arrow](start) -- (step1);
\draw [arrow](step1) -- (step2);
\draw [arrow](step2) -- (step3);
\draw [arrow](step3) -- (step4);
\draw [arrow](step4) -- (step5);
\draw [arrow](step5) -- node[anchor=east] {yes} (step6);
\draw [larrow] (step5.east) -- ++(2,0) node[anchor=south, xshift = -30pt]{no} node
               (lowerright1){} |- (step3.east);
\draw [arrow](step6) -- (step7);
\draw [line width=1pt] (step5.west) -- ++(-1.7,0)
      node(lowerleft1){} |- node[anchor=east, yshift = -35pt]{\fbox{{\color{orange}\shortstack[c]
      {i\\n\\n\\e\\r \\ \\ \\ i\\t\\e\\r}}}} (step3.west);
\draw [arrow](step7) -- node[anchor=east] {yes} (step8);
\draw [larrow] (step7.east) -- ++(2.5,0) node[anchor=south, xshift = -50pt]{no}
      node(lowerright2){} |- (step1.east);

\draw [line width=1pt] (step6.west) -- ++(-3.7,0)
      node(lowerleft1){} |- node[anchor=east, yshift = -90pt]{\fbox{{\color{red}\shortstack[c]
      {o\\u\\t\\e\\r \\ \\ \\ i\\t\\e\\r}}}} (step1.west);
\end{tikzpicture}
\end{center}
\end{figure}

We see from \textbf{Algorithm} \ref{Algorithm B} that there are two level iterations here: the outer iteration is for mesh
refinement, and the inner loop for optimization along each orbital in parallel.
Various optimization approaches can be applied to solve the optimization subproblems, such as
the Newton's method,  the steepest
descent method, and the conjugate gradient method, etc..
In this paper, we focus on the first order algorithms. In our experiments, the choices of search
direction and step size shown  are based on the steepest descent method.

\subsection{Search direction and step size}

We see that the search direction and step size are the most important parts
of an optimization algorithm, which are critical to the speed of convergence. Here, we choose the approximated gradient of the energy function on the Stiefel
manifold as the search direction, and step size is determined by the Barzilai-Borwein
(BB) formula \cite{BaBo}.

Note that two gradient type optimization algorithms are proposed in \cite{ZZWZ}. One of them is OptM-QR
that
calculates $\widetilde{\mathcal{W}} = \mathcal{W} - \tau\nabla E(\mathcal{W})$, and then
orthogonalize $\widetilde{\mathcal{W}}$ by using the QR factorization, i.e.
$\mathcal{W}^{\text{new}} = \text{qr}(\mathcal{W})$.

We see that if the orbitals $\{w_i\}_{i=1}^{N}$ are the eigenfuctions of the Hamiltonian
$H(\rho_{\mathcal{W}})$, the matrix $ \Sigma := \langle (E'(\mathcal{W}))^T\mathcal{W}\rangle $
will be diagonal. Based on this observation, we choose the diagonal element
$\sigma_{ii}$ of $\Sigma$ to approximate the matrix, that is we choose $z_i^{(l)} =
H(\mathcal{W}^{(l)})w_i^{(l)} - \sigma_{ii}^{(l)}w_i^{(l)}$ as our search direction,
where $\sigma_{ii}^{(l)} = \langle H(\mathcal{W}^{(l)})w_i^{(l)},w_i^{(l)}\rangle$.
Then we can calculate the search direction for each orbital in parallel.

In step $l$, we choose the following BB step size (we use the same step size for all
orbitals):
\begin{equation*}
\tau_i^{l,1} = \frac{\text{tr}\langle(\mathcal{S}^{(l-1)})^T\mathcal{S}^{(l-1)}\rangle}
{\text{tr}|\langle(\mathcal{S}^{(l-1)})^T\mathcal{Y}^{(l-1)}\rangle|}
\ \ \text{or} \ \ \tau_i^{l,2} = \frac{\text{tr}|\langle(\mathcal{S}^{(l-1)})^T
\mathcal{Y}^{(l-1)}\rangle|}{\text{tr}\langle(\mathcal{Y}^{(l-1)})^T\mathcal{Y}^{(l-1)}\rangle},
\end{equation*}
where $\mathcal{S}^{(l-1)} = \mathcal{W}^{(l)} - \mathcal{W}^{(l-1)}$, $\mathcal{Y}^{(l-1)} =
\mathcal{Z}^{(l)} - \mathcal{Z}^{(l-1)}$. The calculations of the diagonal elements above can be
done in parallel for each orbital.

In order to guarantee convergence, we perform back tracing for the step size. Let
$C^{(0)} = E(\mathcal{W}^{(0)})$, $Q^{(l+1)} = \eta Q^{(l)} + 1$, $Q^{(0)} = 1$. Then
$\mathcal{W}^{(l+1)} = \mathcal{W}(\tau^{(l)})$, where the $\mathcal{W}(\tau^{(l)})$ is the
new point we get in \textbf{Algorithm} \ref{Algorithm B} in step 4-5 ($\tau^{(l)}_i = \tau^{(l)}$),
and $\tau^{(l)} = \tau_i^{l,1}\delta^s$ or $\tau^{(l)} = \tau_i^{l,2}\delta^s$,
$s$ is smallest nonnegative integer satisfying
\begin{equation}\label{ene-de}
E(\mathcal{W}(\tau^{(l)})) \leq C^{(l)} - \rho_1\tau^{(l)}\|\mathcal{Z}^{(l)}\|_F^2,
\end{equation}
where $C^{(l+1)} = (\eta Q^{(l)}C^{(l)} + E(\mathcal{W}^{(l+1)}))/Q^{(l+1)}$, $\eta$,
$\delta$ and $\rho_1$ are constant factors, we choose $\rho_1 = 1\times10^{-4}$,
$\delta = 0.1$ and $\eta = 0.85$ in our numerical experiments.

\textbf{Algorithm} \ref{Algorithm B} with the above implementation details is named as \textbf{Algorithm} \ref{Algorithm E}.

\begin{algorithm}
\caption{}
\label{Algorithm E}
\begin{algorithmic}[1]
\STATE Given initial space $V^{(0)}$, initial data $\mathcal{U}^{(0)}, \ s.t. \
       \langle(\mathcal{U}^{(0)})^T \mathcal{U}^{(0)}\rangle = I_N$, let $n = 0$;
\STATE Apply an adaptive approach to $ \mathcal{U}^{(n)}$ and get $V^{(n+1)}$;
\STATE $\mathcal{W}^{(0)} = \mathcal{U}^{(n)}$, $l = 0$;
\STATE Calculate the search directions $\{z_i^{(l)}\}_{i=1}^N$ in parallel;
\STATE Calculate the step size $\tau_i^{l} = \tau_i^{l,1}\delta^s$ or $\tau_i^{l} = \tau_i^{l,2}
       \delta^s$, where $s$ is the smallest nonnegative integer such that \eqref{ene-de} is satisfied;
\STATE $\mathcal{W}^{(l+1)} = \text{Orth}(\widetilde{\mathcal{W}}^{(l+1)}), \ l = l + 1$;
\STATE Convergence check for the inner iteration: if not converged, goto step 4;
\STATE $\mathcal{U}^{(n+1)} = \mathcal{W}^{(l)}$, $n = n + 1$;
\STATE Convergence check for the outer iteration: if not converged, goto step 2;
       else, stop.
\end{algorithmic}
\end{algorithm}

\begin{figure}
\centering
\includegraphics[width=0.45\textwidth]{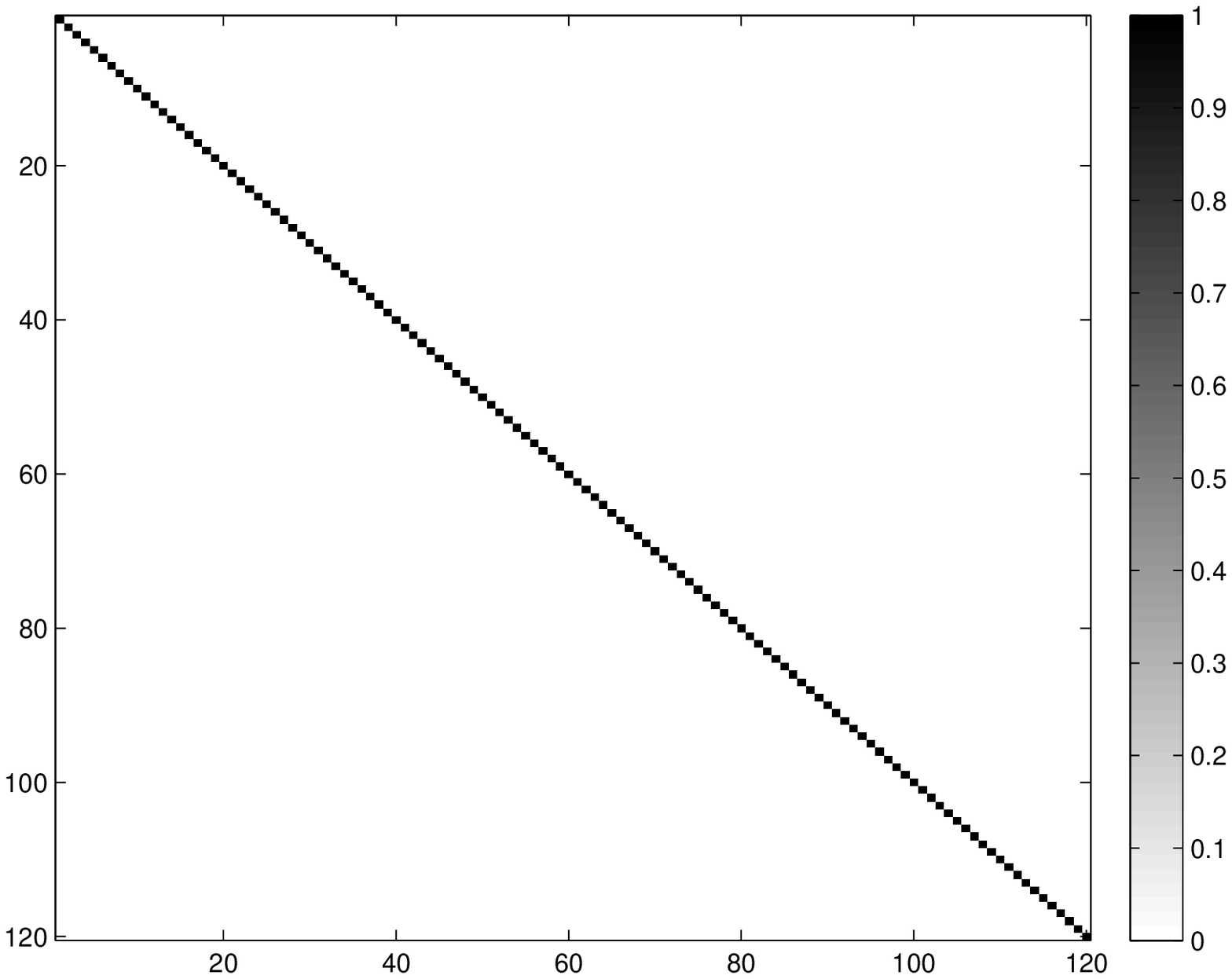}
\includegraphics[width=0.45\textwidth]{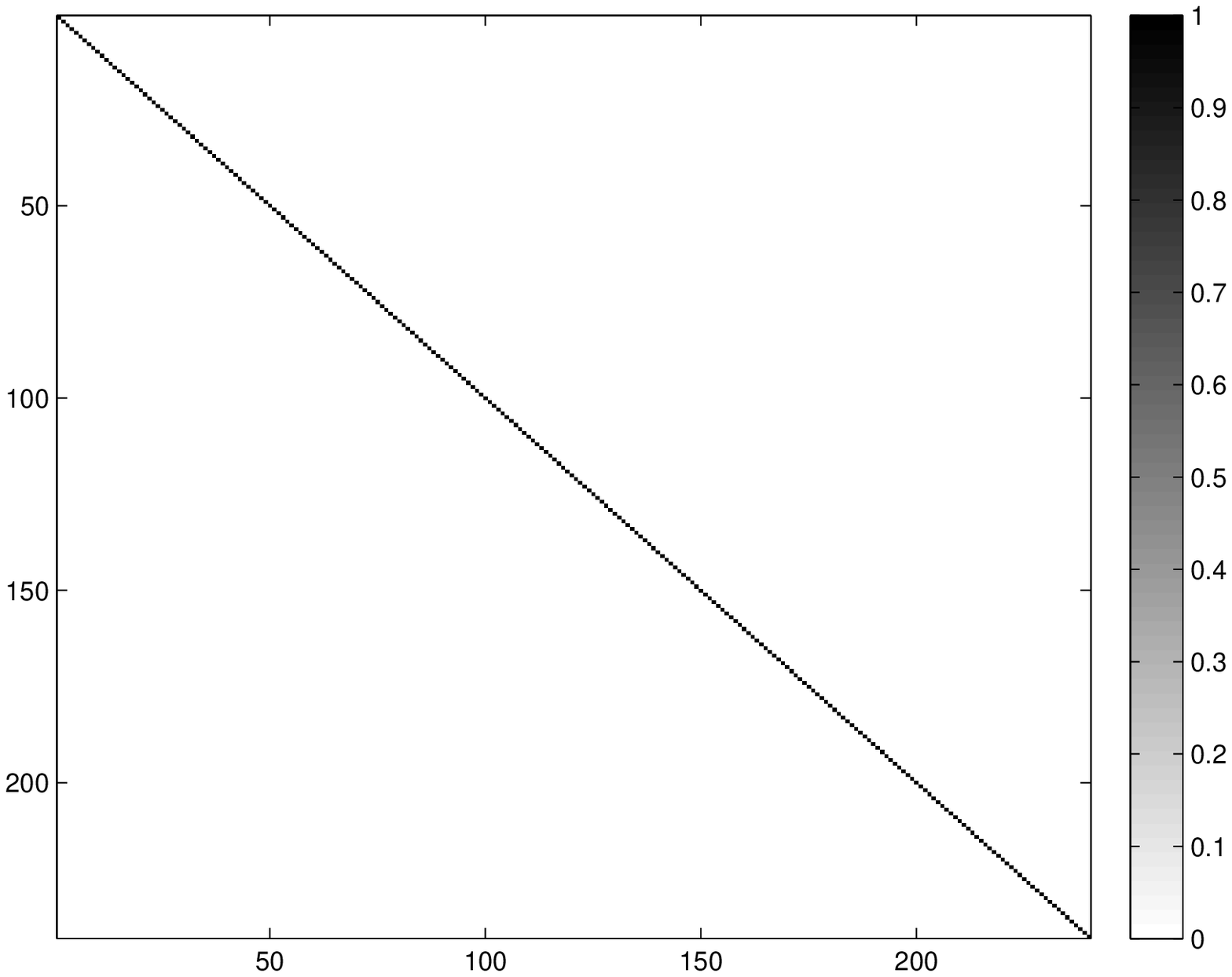} \\
\includegraphics[width=0.45\textwidth]{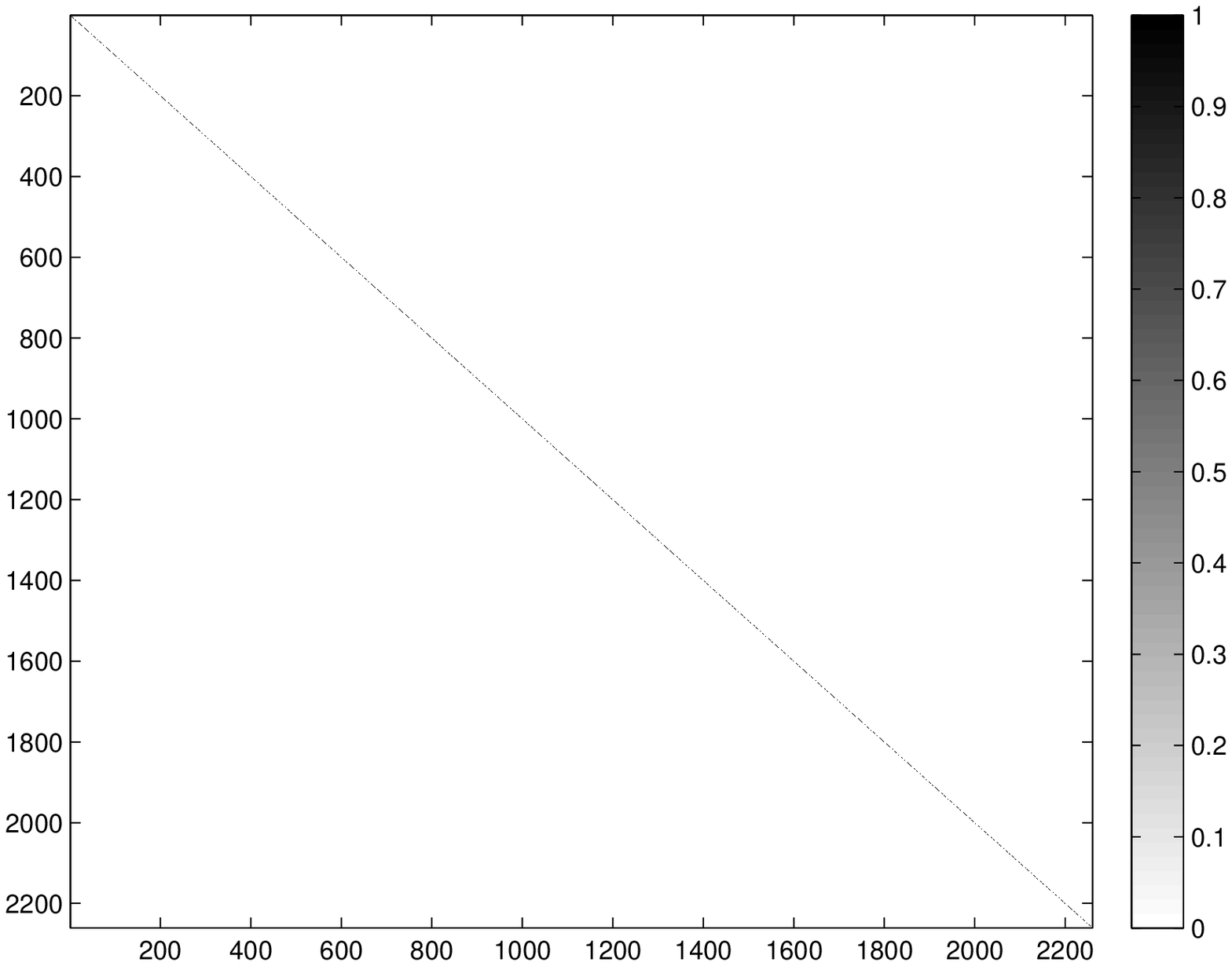}
\includegraphics[width=0.45\textwidth]{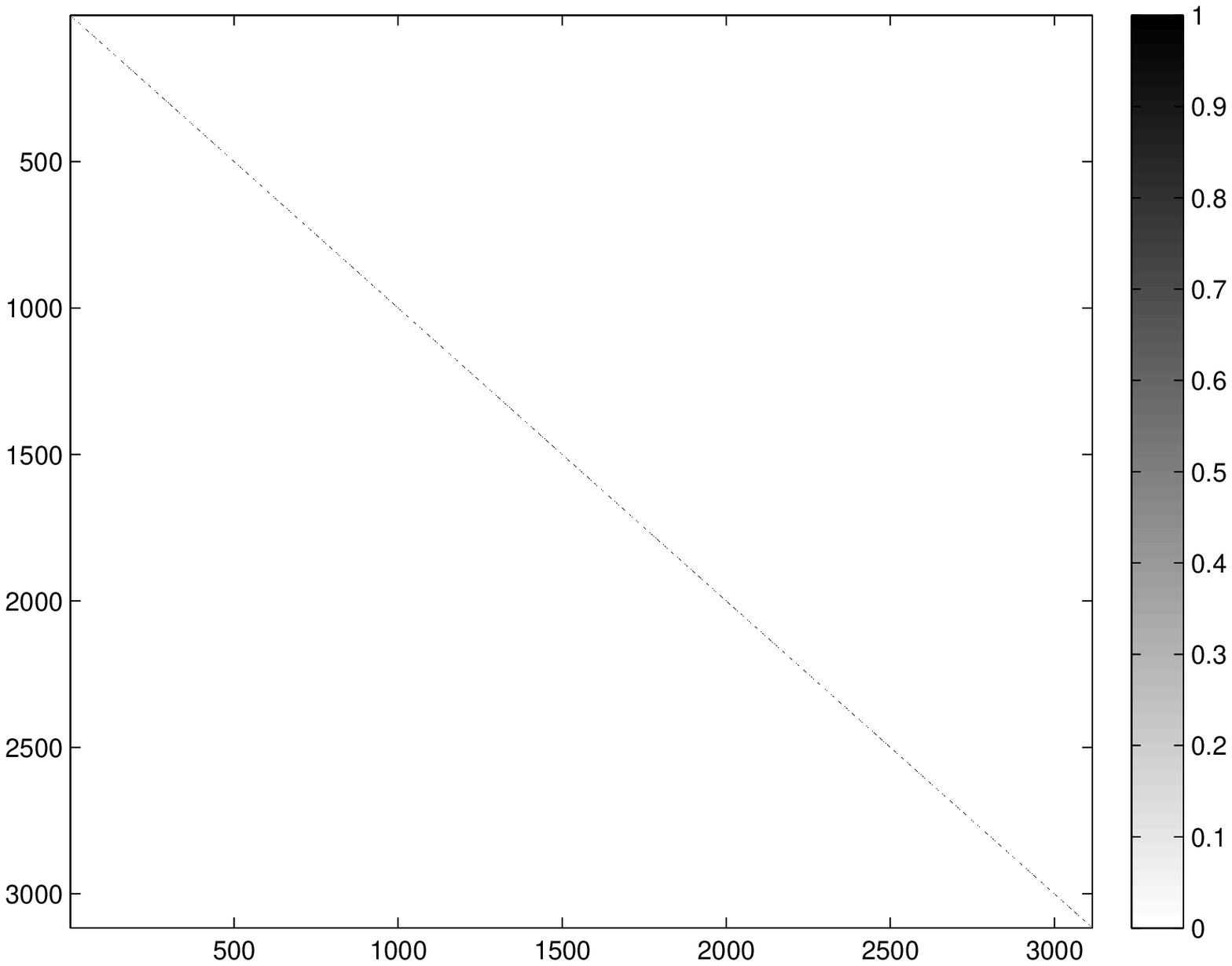}
\caption{The matrix $\langle\widetilde{\mathcal{W}}^T
\widetilde{\mathcal{W}}\rangle$ for $C_{60}$ and $C_{120}$, $C_{1015}H_{460}$ and \ $C_{1419}H_{556}$
by \textbf{Algorithm} \ref{Algorithm M}, at iteration step 101, with $n_{diag}$ = 100, $n_{orthg}$ = 2.}
\label{finorm}
\end{figure}

As in Proposition 3.2 in \cite{ZZWZ}, if we use the full $\Sigma$ for all the orbitals, the
resulting $\widetilde{\mathcal{W}}$ is full rank and $\langle\widetilde{\mathcal{W}}^T
\widetilde{\mathcal{W}}\rangle$ is well conditioned, which makes it easy to perform QR factorization
based on Cholesky factorization. While for the algorithm proposed here, we do not have the theoretical
result to ensure this property. However, this is not a problem to perform the QR factorization
numerically (if we use the subspace diagonalization technique introduced in the next subsection),
and  we observe from Figure \ref{finorm} that the matrix $\langle\widetilde{\mathcal{W}}^T
\widetilde{\mathcal{W}}\rangle$ is diagonal dominated and very close to identity matrix.

\subsection{Modifications}

As the iteration goes on, the matrix $\Sigma := \langle (E'(\mathcal{W}))^T\mathcal{W}\rangle$
deviates from diagonal matrix gradually, the error to approximate the gradient by choosing the
diagonal elements of $\Sigma$ becomes larger and larger. Therefore, we perform subspace diagonalization
for the Hamiltonian in the space spanned by current orbitals every $n_{diag}$ (a parameter) steps.
Specifically, suppose orthogonal matrix $P\in \mathbb{R}^{N\times N}$ diagonalize $\Sigma$, i.e.
$P^T\Sigma P = \Gamma$ ($\Gamma$ is diagonal), let $\mathcal{W}^{\text{new}} = \mathcal{W}P$. Then if
$\langle\mathcal{W}^T\mathcal{W}\rangle = I_N$, we have $\langle(\mathcal{W}^{\text{new}})^T
\mathcal{W}^{\text{new}}\rangle = P^T\langle(\mathcal{W}^T\mathcal{W})\rangle P = I_N$,
$\rho(\mathcal{W}^{\text{new}}) = \rho(\mathcal{W})$. Since $H(\mathcal{W})$ is determined by
$\rho(\mathcal{W})$, we have $H(\mathcal{W}^{\text{new}}) = H(\mathcal{W})$, then $\Sigma^{\text{new}}
:= \langle(H(\mathcal{W}^{\text{new}})\mathcal{W}^{\text{new}})^T \mathcal{W}^{\text{new}}\rangle$
is diagonal.

We see from \textbf{Algorithm} \ref{Algorithm E} that the main part of calculations we cannot carry out
in parallel for each orbital are the orthogonalization of orbitals in step 6. To decrease this
part of time, we perform orthogonalization every $n_{org}$ steps. That is, we do not orthogonalize
the orbitals every iteration. As a result, we calculate the total energy and perform back tracing
for the step size every $n_{org}$ steps. Such a modified algorithm is stated as \textbf{Algorithm}
\ref{Algorithm M}:

\begin{algorithm}
\caption{}
\label{Algorithm M}
\begin{algorithmic}[1]
\STATE Given initial space $V^{(0)}$, initial data $\mathcal{U}^{(0)}, \ s.t. \ \langle(\mathcal{U}^{(0)})^T
       \mathcal{U}^{(0)}\rangle = I_N$, $n_{orthg}$ and $n_{diag}$, let $n = 0$;
\STATE Apply an adaptive approach to $\mathcal{U}^{(n)}$ and get $V^{(n+1)}$;
\STATE $\mathcal{W}^{(0)} = \mathcal{U}^{(n)}$, $l = 0$;
\STATE Calculate the search directions $\{z_i^{(l)}\}_{i=1}^N$ in parallel;
\IF {(mod(l, $n_{org}$) $\ne$ 0)}
\STATE Calculate the step size $\tau_i^{l} = \tau_i^{l,1}$ or $\tau_i^{l} = \tau_i^{l,2}$, let
       $w_i^{(l+1)} = w_i^{(l)} - \tau_i^{(l)}z_i^{(l)}$, $l = l + 1$;
\ELSE
\STATE Calculate the step size $\tau_i^{l} = \tau_i^{l,1}\delta^s$ or $\tau_i^{l} = \tau_i^{l,2}
       \delta^s$, where $s$ is the smallest nonnegative integer such that \eqref{ene-de} is satisfied;
\STATE $\mathcal{W}^{(l+1)} = \text{Orth}(\widetilde{\mathcal{W}}^{(l+1)}), \ l = l + 1$;
\STATE If $mod(l,n_{diag}) = 0$, let $\mathcal{W}^{(l+1)} = \mathcal{W}^{(l+1)}P^{(l+1)}$,
       where $P^{(l+1)}$ is the rotation matrix introduced in the previous subsection;
\ENDIF
\STATE Convergence check for the inner iteration: if not converged, goto step 4;
\STATE $\mathcal{U}^{(n+1)} = \mathcal{W}^{(l)}$, $n = n + 1$;
\STATE Convergence check for the outer iteration: if not converged, goto step 2;
       else, stop.
\end{algorithmic}
\end{algorithm}

\section{Numerical experiments}

Our numerical experiments are carried out on the software package OCTOPUS\footnote{www.tddft.org/programs/
octopus.}(version 4.0.1), a real space ab initio computing platform. We choose local density
approximation (LDA) to approximate $v_{xc}(\rho)$ \cite{PeZu} in our numerical experiments.
The initial guess of the orbitals is generated by linear combination of the atomic orbits (LCAO)
method, and the Troullier-Martins norm conserving pseudopotential \cite{TrMa} is used
in our computation. We use \textbf{Algorithm \ref{Algorithm M}} for all the tests, and the current
results are without mesh refinement.

Our examples include several typical molecular systems: benzene ($C_6H_6$), aspirin ($C_9H_8O_4)$,
fullerene ($C_{60}$), alanine chain $(C_{33}H_{11}O_{11}N_{11})$, carbon
nano-tube ($C_{120}$), biological ligase 2JMO $(C_{178}H_{283}O_{50}N_{57}S)$,
protein fasciculin2 $(C_{276}H_{442}O_{90}N_{88}S_{10})$, carbon cluster $C_{1015}H_{460}$ and
$C_{1419}H_{556}$. 

\subsection{Reliability}
We first test the reliability of out algorithm. The results are shown in Table \ref{t1}.
The Opt-Z3W is the OptM-QR algorithm in \cite{ZZWZ}, Opt-Par is the parallel orbital-updating
based optimization algorithm proposed in this paper. The column 'energy' is the Kohn-Sham
total energy (in atomic units), 'iter' denotes the total number of iterations, and 'resi'
denotes the norm of the residual $\|\nabla E(\mathcal{U})\| = \sqrt{\text{tr}\langle
\nabla E(\mathcal{U})^T\nabla E(\mathcal{U})\rangle}$. We use 'cores' to denote the number of
CPU cores used in that calculation, and it is chosen as $2^s$, where $s$ is the largest integer
that makes the number of elements on each processor not smaller than the value recommended by
Octopus.

We choose $n_{diag} = 100$, $n_{org} = 1$. Since the mesh refinement is not implemented,
we use $\|\nabla E(\mathcal{U})\| < 10^{-6}$ as the convergence criterion for the inner iteration.
And since it is difficult to converge to this criterion for large systems, we choose a looser criterion,
that is: the average of absolute value of each element of the matrix $\nabla E(\mathcal{U})$ is less
than $5\times10^{-9}$. So, for 2JMO, FAS2, $C_{1015}H_{460}$ and $C_{1419}H_{556}$, the practical
convergence criteria are $4.06\times10^{-5}, 6.48\times10^{-5}, 7.62\times10^{-5} \ \mbox{and} \
1.01\times10^{-4}$.

\begin{center}
\begin{table}[!htbp]
\caption{The results using $\|\nabla E\|_F$ as
         convergence criteria, $n_{diag} = 100$, $n_{org} = 1$.}\label{t1}
\begin{center}
\begin{tabular}{|c| c c c c|}
\hline
algorithm & energy & iter & resi & total-time(s)\\
\hline
\multicolumn{5}{|c|}{benzene($C_6H_6) \ \ \  N_g = 102705 \ \ \  N = 15 \ \ \  cores = 8$} \\
\hline
Opt-Z3W & -3.74246025E+01 & 164 & 7.49E-07 & 6.63   \\
Opt-Par & -3.74246025E+01 & 173 & 7.95E-07 & 7.58   \\
\hline
\multicolumn{5}{|c|}{aspirin($C_9H_8O_4) \ \ \  N_g = 133828 \ \ \  N = 34 \ \ \  cores = 16$} \\
\hline
Opt-Z3W     & -1.20214764E+02 & 153 & 9.89E-07 & 15.67 \\
Opt-Par     & -1.20214764E+02 & 134 & 9.55E-07 & 17.33 \\
\hline
\multicolumn{5}{|c|}{$C_{60} \ \ \ N_g = 191805  \ \ \  N = 120 \ \ \  cores = 16$} \\
\hline
Opt-Z3W     & -3.42875137E+02 & 222   & 9.02E-07 & 96.73  \\
Opt-Par     & -3.42875137E+02 & 230   & 9.43E-07 & 108.24 \\
\hline
\multicolumn{5}{|c|}{alanine chain$(C_{33}H_{11}O_{11}N_{11})\ \ \  N_g = 293725 \ \ \  N = 132 \ \ \  cores = 32$} \\
\hline
Opt-Z3W & -4.78562217E+02 & 1779 & 9.60E-07 & 1036.44 \\
Opt-Par & -4.78562217E+02 & 1673 & 9.87E-07 & 1101.63 \\
\hline
\multicolumn{5}{|c|}{$C_{120} \ \ \ N_g = 354093 \ \ \  N = 240 \ \ \  cores = 32$} \\
\hline
Opt-Z3W & -6.84467048E+02 & 2392 & 9.64E-07 & 2560.63 \\
Opt-Par & -6.84467048E+02 & 1848 & 9.42E-07 & 2274.45 \\
\hline
\multicolumn{5}{|c|}{2JMO$(C_{178}H_{283}O_{50}N_{57}S) \ \ \  N_g = 1226485 \ \ \  N = 793 \ \ \  cores = 128$} \\
\hline
Opt-Z3W & -2.56413551E+03 & 1502  & 4.03E-05 & 11998.30  \\
Opt-Par & -2.56413551E+03 & 1846  & 4.05E-05 & 18201.79  \\
\hline
\multicolumn{5}{|c|}{FAS2$(C_{276}H_{442}O_{90}N_{88}S_{10}) \ \ \  N_g = 1903841 \ \ \ N = 1293 \ \ \  cores = 256$} \\
\hline
Opt-Z3W  & -4.26018875E+03 & 2050  & 6.46E-05 & 30152.76   \\
Opt-Par  & -4.26018881E+03 & 2094  & 6.48E-05 & 44952.46  \\
\hline
\multicolumn{5}{|c|}{$C_{1015}H_{460}\ \ \  N_g = 1462257 \ \ \ N = 2260 \ \ \  cores = 256$} \\
\hline
Opt-Z3W  & -6.06369982E+03 & 137   & 3.08E-05 & 5331.19   \\
Opt-Par  & -6.06369982E+03 & 144   & 7.26E-05 & 6217.72  \\
\hline
\multicolumn{5}{|c|}{$C_{1419}H_{556}\ \ \  N_g = 1828847 \ \ \ N = 3116 \ \ \  cores = 512$} \\
\hline
Opt-Z3W  & -8.43085432E+03 & 151   & 1.01E-04 & 9084.88   \\
Opt-Par  & -8.43085432E+03 & 138   & 5.49E-05 & 9453.76  \\
\hline
\end{tabular}
\end{center}
\end{table}
\end{center}

We see from Table \ref{t1} that for all the systems, Opt-Par get the same convergence
result as Opt-Z3W, which confirms the reliability of our algorithm.

\begin{figure}
\centering
\includegraphics[width=0.45\textwidth]{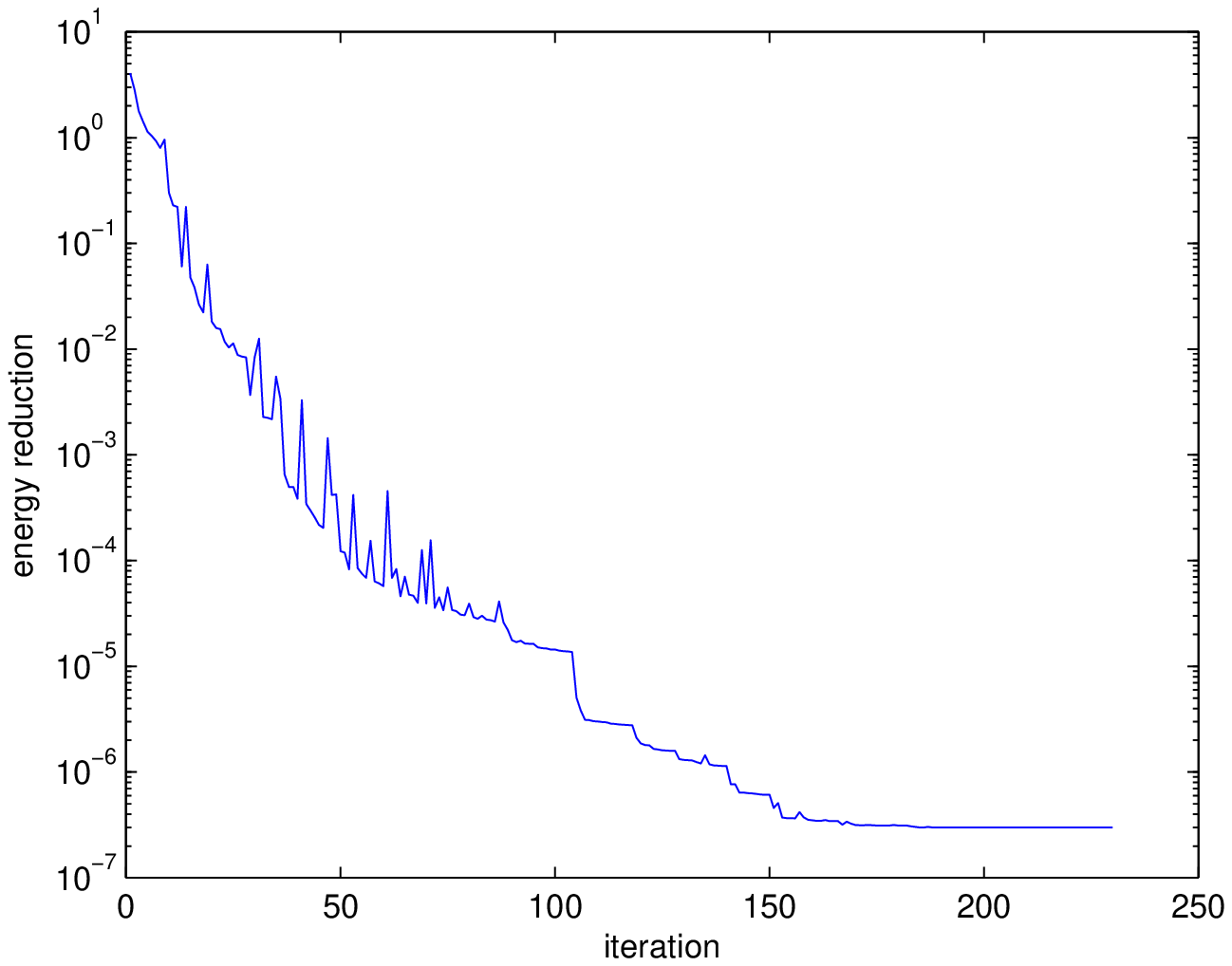}
\includegraphics[width=0.45\textwidth]{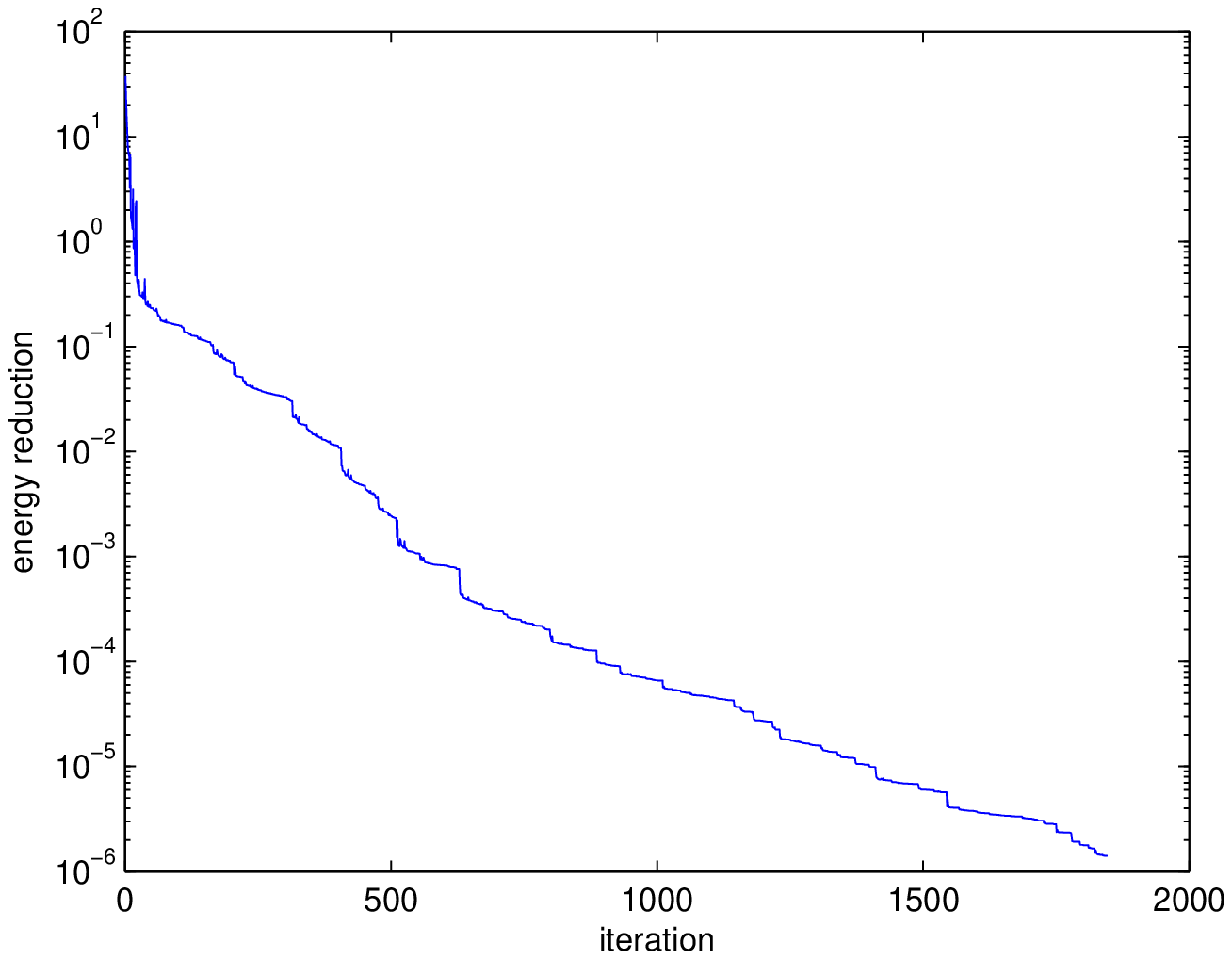}
\caption{Energy reduction for C60 and 2JMO.}
\label{f1}
\end{figure}

We  show the energy reduction $E(X^{(n)}) - E_{min}$ for $C_{60}$ and 2JMO in Figure
\ref{f1}, where $E_{min}$ is a reliable minimum of the total energy. We can see
that although there are some oscillations, the energy decreases gradually until convergence
is reached.

\subsection{Parallelism}

In this subsection, we test the parallelism in orbital for our algorithm. Since the parallel
in orbital is not implemented yet, we show the computation time in two parts: the first part
is the total computation time ('total-time'), the second part is the time that can be done
parallel in each orbital ('par-time'), 'percent' is 'par-time' divided by 'total-time'.
Because it costs too much time to calculate the residual $\|\nabla E(\mathcal{U})\|$, and this
part of calculations are not parallel in orbital, we use the relative
change of total energy as our convergence criteria, that is: $\frac{df^{(l)}+df^{(l-1)}+
df^{(l-2)}}{3}<10^{-13}$, where $df^{(l)} = \frac{|E(\mathcal{U}^{(l-1)}) - E(\mathcal{U}^{(l)})|}
{|E(\mathcal{U}^{(l-1)})|+1}$ (when we perform orthogonalization every $n_{org}$ steps, since we don't
calculate the energy step, we choose the average change of total energy of adjacent $n_{org}$ steps).
We present the results under different parameters in the Tables \ref{t3} and \ref{t4}.

\begin{center}
\begin{table}[!htbp]
\caption{The results using convergence criteria $\frac{df^{(n)}+df^{(n-1)}+df^{(n-2)}}
        {3}<10^{-13}$, $n_{diag} = 100$, $n_{org} = 1$.}\label{t3}
\begin{center}
\begin{tabular}{|>{\footnotesize}c|>{\small}c >{\small}c >{\small}c >{\small}c >{\small}c >{\footnotesize}c|}
\hline
algorithm & energy & iter & resi & total-time(s) & par-time(s) & percent \\
\hline
\multicolumn{7}{|c|}{benzene$ \ \ \  N_g = 102705 \ \ \  N = 15 \ \ \  cores = 8$} \\
\hline
 Opt-Z3W & -3.74246025E+01 & 165 & 7.93E-07 & 6.63   & - & - \\
 Opt-Par & -3.74246025E+01 & 168 & 2.87E-06 & 6.94   & 3.86  & 56\%\\
\hline
\multicolumn{7}{|c|}{aspirin$ \ \ \  N_g = 133828 \ \ \  N = 34 \ \ \  cores = 16$} \\
\hline
 Opt-Z3W & -1.20214764E+02 & 141 & 1.96E-06 & 15.67   & - & - \\
 Opt-Par & -1.20214764E+02 & 135 & 2.96E-06 & 14.52 & 8.23 & 57\%\\
\hline
\multicolumn{7}{|c|}{$C_{60} \ \ \ N_g = 191805  \ \ \  N = 120 \ \ \  cores = 16$} \\
\hline
 Opt-Z3W & -3.42875137E+02 & 176 & 7.69E-06 & 99.66 & - & - \\
 Opt-Par & -3.42875137E+02 & 191 & 4.00E-05 & 80.64 & 54.44  & 68\% \\
\hline
\multicolumn{7}{|c|}{alanine chain$\ \ \  N_g = 293725 \ \ \  N = 132 \ \ \  cores = 32$} \\
\hline
 Opt-Z3W & -4.78562217E+02 & 1538 & 1.73E-05 & 968.12 & - & - \\
 Opt-Par & -4.78562217E+02 & 1599 & 4.48E-05 & 942.87 & 576.32 & 61\%\\
\hline
\multicolumn{7}{|c|}{$C_{120} \ \ \ N_g = 354093 \ \ \  N = 240 \ \ \  cores = 32$} \\
\hline
 Opt-Z3W & -6.84467048E+02 & 1390 & 1.80E-05 & 1588.33 & - & - \\
 Opt-Par & -6.84467048E+02 & 1435 & 4.53E-05 & 1403.23 & 1056.02 & 75\%\\
\hline
\multicolumn{7}{|c|}{2JMO$ \ \ \  N_g = 1226485 \ \ \  N = 793 \ \ \  cores = 128$} \\
\hline
 Opt-Z3W & -2.56413551E+03 & 1480  & 4.08E-05 & 14524.26  & - & - \\
 Opt-Par & -2.56413551E+03 & 1794  & 4.18E-05 & 13692.19  & 9564.53 & 70\%\\
\hline
\multicolumn{7}{|c|}{FAS2$ \ \ \  N_g = 1903841 \ \ \ N = 1293 \ \ \  cores = 256$} \\
\hline
 Opt-Z3W & -4.26018875E+03 & 2150  & 5.67E-05 & 43218.37  & - & - \\
 Opt-Par & -4.26018881E+03 & 2225  & 5.83E-05 & 37925.49  & 25147.62 & 66\%\\
\hline
\multicolumn{7}{|c|}{$C_{1015}H_{460}\ \ \  N_g = 1462257 \ \ \ N = 2260 \ \ \  cores = 256$} \\
\hline
Opt-Z3W  & -6.06369982E+03 & 145   & 9.59E-06 & 5525.22  & - & - \\
Opt-Par  & -6.06369982E+03 & 157   & 1.41E-04 & 4873.35  & 3012.72 & 62\% \\
\hline
\multicolumn{7}{|c|}{$C_{1419}H_{556}\ \ \  N_g = 1828847 \ \ \ N = 3116 \ \ \  cores = 512$} \\
\hline
Opt-Z3W  & -8.43085432E+03 & 165   & 5.05E-05 & 9926.59  & - & - \\
Opt-Par  & -8.43085432E+03 & 179   & 6.43E-05 & 9277.48  & 6223.67 & 67\%\\
\hline
\end{tabular}
\end{center}
\end{table}
\end{center}

\begin{figure}
\centering
\includegraphics[width=0.45\textwidth]{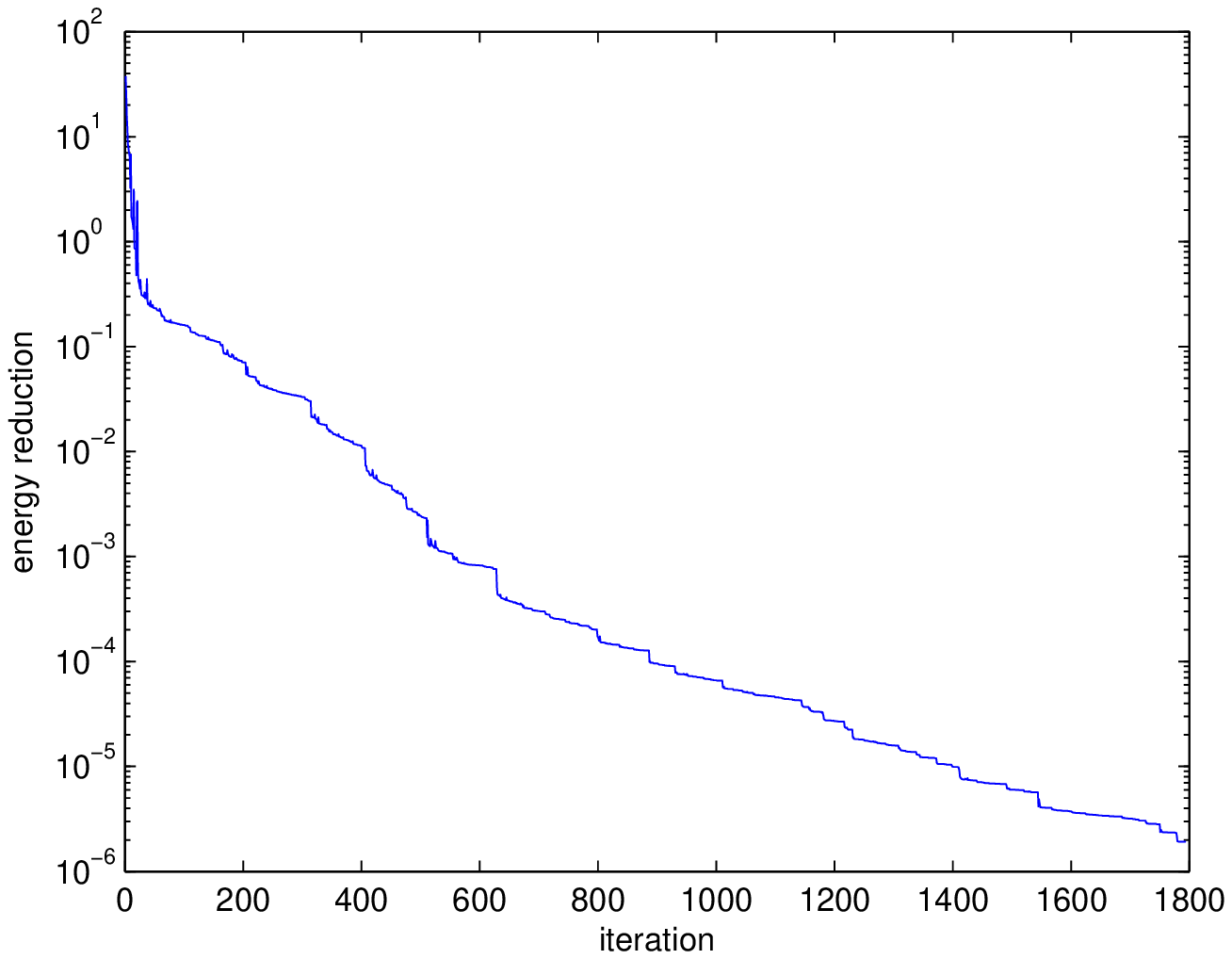}
\includegraphics[width=0.45\textwidth]{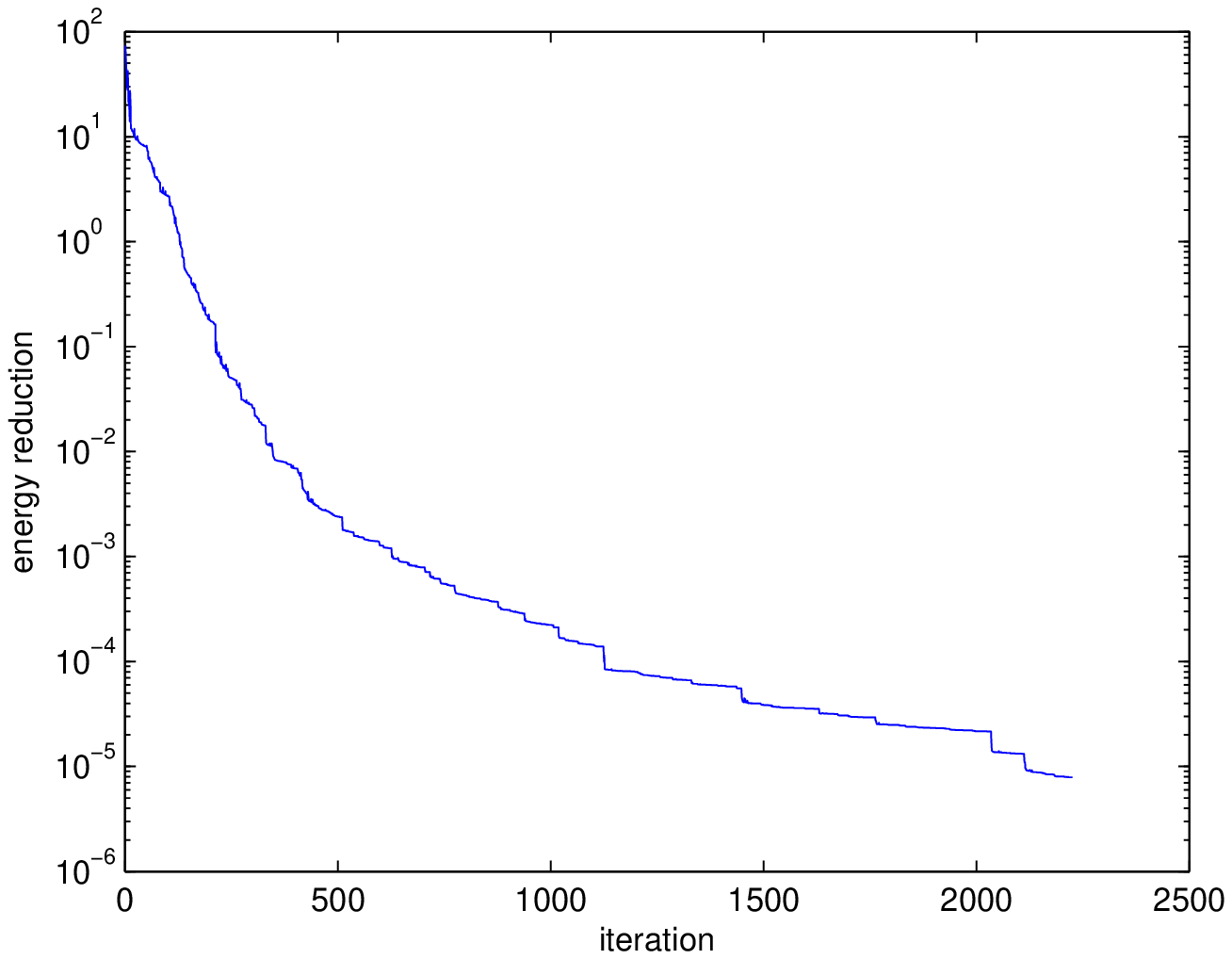}
\caption{Energy reduction for 2JMO and FAS2 (Correspond to Table \ref{t3}).}
\label{f2}
\end{figure}

\begin{center}
\begin{table}[!htbp]
\caption{The results using convergence criteria
         $\frac{df^{(n)}+df^{(n-1)}+df^{(n-2)}}{3}< 5 \times 10^{-13}$, $n_{diag} = 50$, $n_{org} = 2$.}\label{t4}
\begin{center}
\begin{tabular}{|>{\footnotesize}c|>{\small}c >{\small}c >{\small}c >{\small}c >{\small}c >{\footnotesize}c|}
\hline
algorithm & energy & iter & resi & total-time(s) & par-time(s) & percent \\
\hline
\multicolumn{7}{|c|}{benzene$ \ \ \  N_g = 102705 \ \ \  N = 15 \ \ \  cores = 8$} \\
\hline
 Opt-Z3W & -3.74246025E+01 & 142 & 9.76E-06 & 5.85   & - & - \\
 Opt-Par & -3.74246025E+01 & 153 & 4.72E-06 & 4.39   & 2.58  & 59\%\\
\hline
\multicolumn{7}{|c|}{aspirin$ \ \ \  N_g = 133828 \ \ \  N = 34 \ \ \  cores = 16$} \\
\hline
 Opt-Z3W & -1.20214764E+02 & 121 & 1.72E-05 & 12.68   & - & - \\
 Opt-Par & -1.20214764E+02 & 119 & 7.99E-06 & 9.61    & 5.55 & 58\%\\
\hline
\multicolumn{7}{|c|}{$C_{60} \ \ \ N_g = 191805  \ \ \  N = 120 \ \ \  cores = 16$} \\
\hline
 Opt-Z3W & -3.42875137E+02 & 167 & 3.46E-05 & 73.18 & - & - \\
 Opt-Par & -3.42875137E+02 & 189 & 4.66E-05 & 84.50 & 39.37  & 72\% \\
\hline
\multicolumn{7}{|c|}{alanine chain$\ \ \  N_g = 293725 \ \ \  N = 132 \ \ \  cores = 32$} \\
\hline
 Opt-Z3W & -4.78562217E+02 & 1017 & 2.90E-05 & 659.92 & - & - \\
 Opt-Par & -4.78562217E+02 & 1865 & 2.12E-05 & 778.51 & 559.43 & 72\%\\
\hline
\multicolumn{7}{|c|}{$C_{120} \ \ \ N_g = 354093 \ \ \  N = 240 \ \ \  cores = 32$} \\
\hline
 Opt-Z3W & -6.84467046E+02 & 1200 & 5.05E-05 & 1396.19 & - & - \\
 Opt-Par & -6.84467048E+02 & 1965 & 1.88E-05 & 1346.37 & 1100.22 & 82\%\\
\hline
\multicolumn{7}{|c|}{2JMO$ \ \ \  N_g = 1226485 \ \ \  N = 793 \ \ \  cores = 128$} \\
\hline
 Opt-Z3W & -2.56413551E+03 & 1049  & 8.98E-05 & 9836.49   & - & - \\
 Opt-Par & -2.56413551E+03 & 1987  & 4.47E-05 & 9593.32   & 7180.12 & 75\%\\
\hline
\multicolumn{7}{|c|}{FAS2$ \ \ \  N_g = 1903841 \ \ \ N = 1293 \ \ \  cores = 256$} \\
\hline
 Opt-Z3W & -4.26018872E+03 & 1596  & 1.23E-04 & 31647.13  & - & - \\
 Opt-Par & -4.26018870E+03 & 2539  & 1.12E-04 & 26709.60  & 19401.66 & 73\%\\
\hline
\multicolumn{7}{|c|}{$C_{1015}H_{460}\ \ \  N_g = 1462257 \ \ \ N = 2260 \ \ \  cores = 256$} \\
\hline
Opt-Z3W  & -6.06369982E+03 & 128   & 9.90E-05 & 4899.40  & - & - \\
Opt-Par  & -6.06369982E+03 & 139   & 3.23E-05 & 2615.67  & 1611.80 & 62\% \\
\hline
\multicolumn{7}{|c|}{$C_{1419}H_{556}\ \ \  N_g = 1828847 \ \ \ N = 3116 \ \ \  cores = 512$} \\
\hline
Opt-Z3W  & -8.43085432E+03 & 150   & 1.02E-04 & 9135.69  & - & - \\
Opt-Par  & -8.43085432E+03 & 155   & 4.43E-05 & 5125.12  & 3143.89 & 61\%\\
\hline
\end{tabular}
\end{center}
\end{table}
\end{center}

\begin{figure}
\centering
\includegraphics[width=0.45\textwidth]{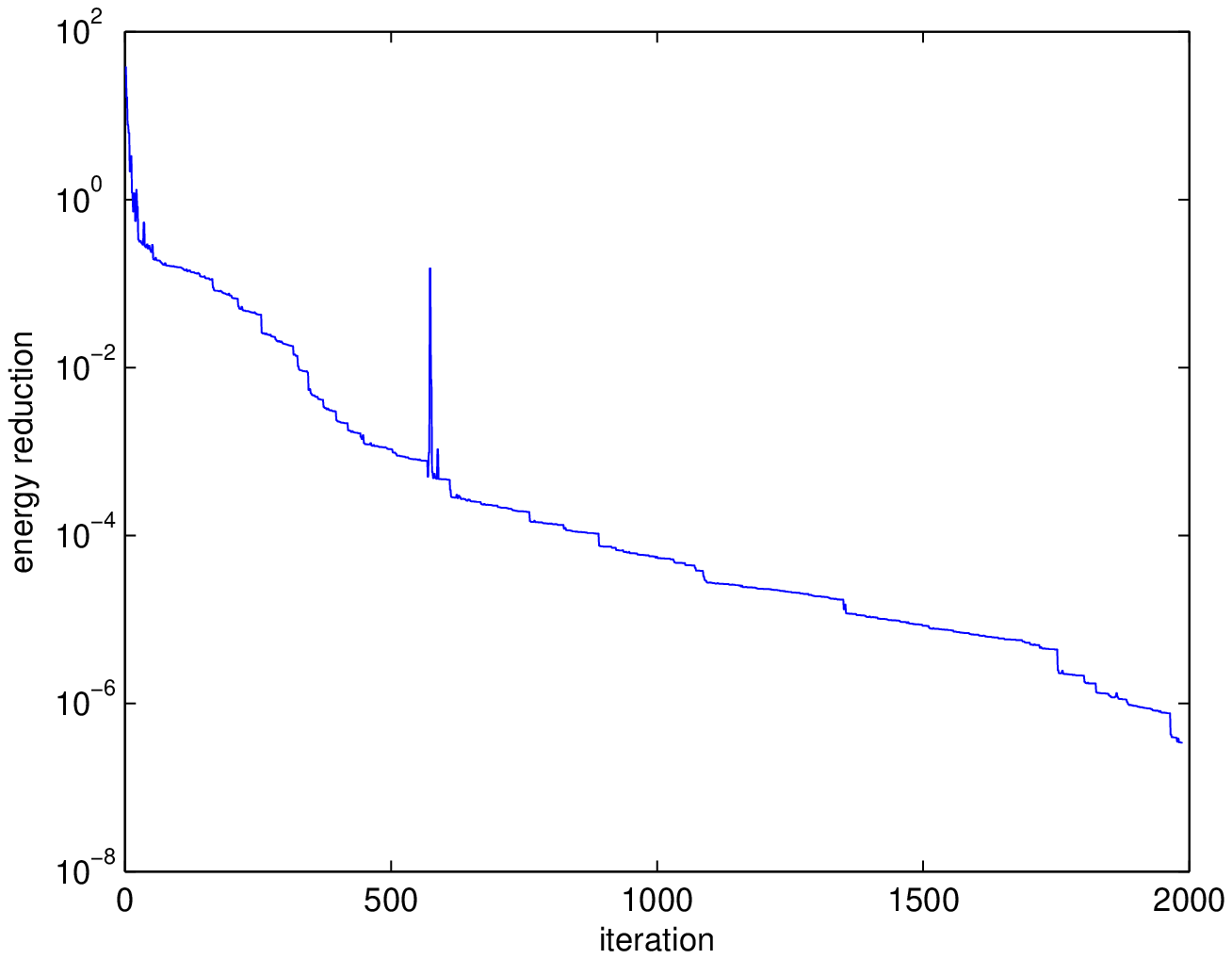}
\includegraphics[width=0.45\textwidth]{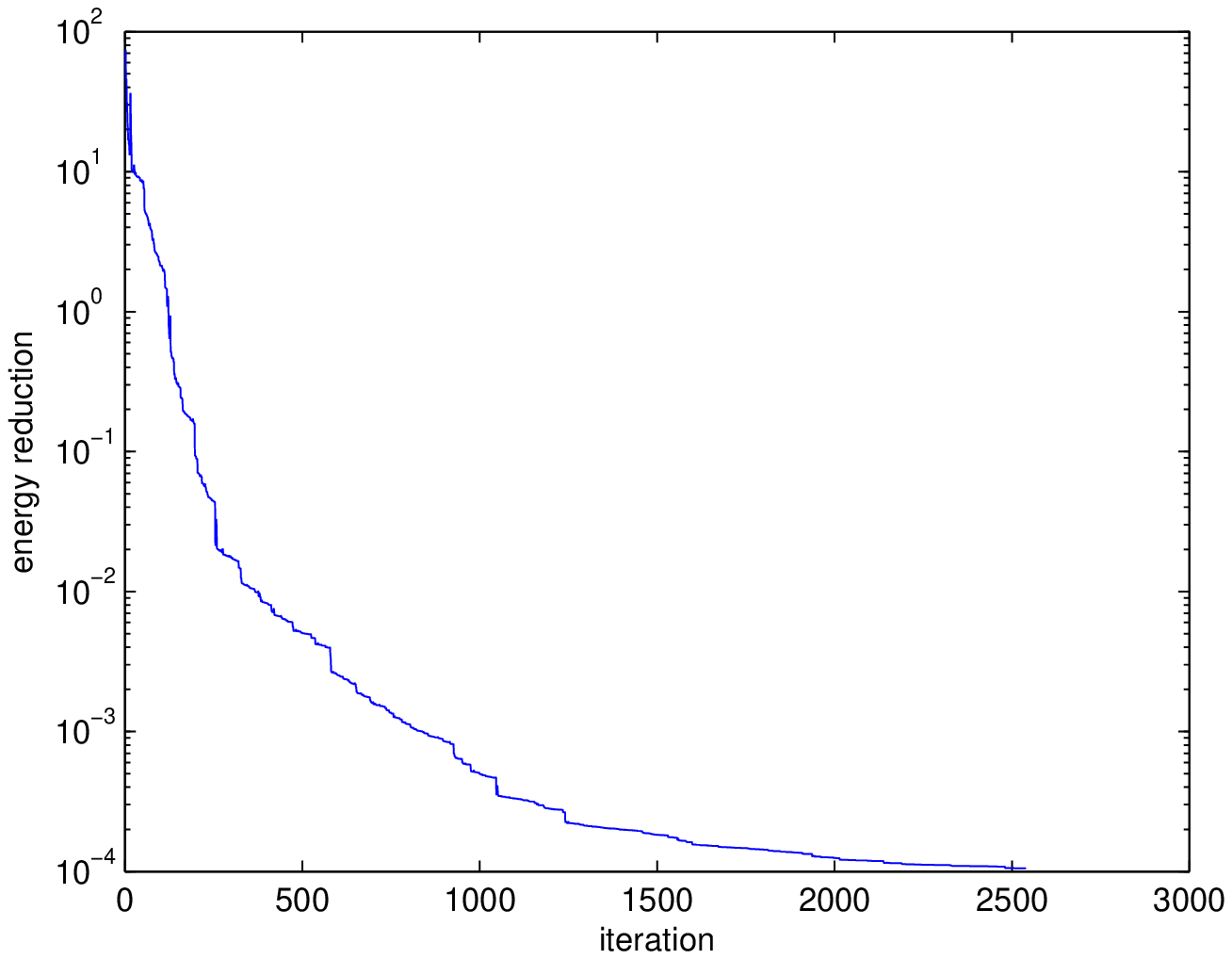}
\caption{Energy reduction for 2JMO and FAS2 (Correspond to Table \ref{t4}).}
\label{f3}
\end{figure}

We observe from Table \ref{t3} that: for larger system, the percentage of time that we can do
parallel in orbital is larger than 65\%, which is very attractive for large systems calculations.
By choosing $n_{org} = 2$, the parallel percentage is further increased, as shown in Table
\ref{t4}.

Figure \ref{f2} and \ref{f3} show the energy reduction for 2JMO and FAS2, which confirm the
convergence of our algorithm. One may see that there are some oscillations in the energy reduction
curve for 2JMO in Figure \ref{f3}. This is due to we do not perform orthogonalization every step,
the search direction may deviate from the energy decrease direction, and we do not perform back
tracing for the step size in these steps. However, the oscillation steps are far form convergence,
and there is no oscillation when the convergence is reached. Therefore, it does not affect the
computational results.

\section{Concluding remarks}
In this paper, we proposed a parallel optimization method for electronic
structure calculations based on a single orbital-updating approximation. It is shown by our numerical experiments  that this new algorithm
is reliable and has very good potential for large scale electronic structure calculations.
We believe that our parallel optimization algorithms should be friendly to supercomputers.
 It is our ongoing work to carry out the two-level parallel version of our algorithms.

\end{document}